\documentclass{article}
\usepackage{amssymb}
\usepackage{amsmath,amsfonts,amsthm}

\newtheorem{algorithm}{Algorithm}

\newtheorem{proposition}{Proposition}
\newtheorem{theorem}{Theorem}
\newtheorem{definition}{Definition}
\newtheorem{remark}{Remark}

\begin{document}

\title{A general framework for applying FGLM techniques to linear codes}

\author{ M. Borges-Quintana \and M. A. Borges-Trenard \thanks{Dpto. de Matem\'atica, FCMC, U. de Oriente, Santiago de Cuba, Cuba. First author is partially supported by Spanish MEC grant SB 2003-0287. \tt{mijail@mbq.uo.edu.cu, mborges@mabt.uo.edu.cu}} \and E. Mart\'{\i}nez-Moro\thanks{Dpto. de Matem\'atica Aplicada, U. de Valladolid, Spain. Partially supported by Spanish MEC MTM2004-00876 and MTM2004-00958. \tt{edgar@maf.uva.es}}}

\maketitle              % typeset the title of the contribution

\begin{abstract}
We show herein that a pattern based on FGLM techniques can be used for computing Gr\"obner bases, or related structures, associated to linear codes. This Gr\"obner bases setting turns out to be strongly related to the combinatorics of the codes.
\end{abstract}
\section*{Introduction}
It is well known that the complexity of Gr\"obner bases computation heavily depends on the term orderings, moreover, elimination
orderings often yield a greater complexity. This remark led to the so
called FGLM convertion problem, i.e., {\bf given} a Gr\"obner basis
with respect to a certain term ordering, {\bf
find} a Gr\"obner basis of the same ideal with respect to another
term ordering. One of the efficient approaches for solving this
problem, in the zero-dimensional case, is the FGLM algorithm (see
\cite{FGLM}).\\ 
The key ideas of this algorithm were successfully generalized in \cite{MMM} with the objective of computing Gr\"obner bases of
zero-dimensional ideals that are determined by functionals. 
%(in the
%sense that they are kernels of finite sets of morphisms from
%the polynomial ring to the base field). 
In fact, the pioneer work of FGLM and \cite{MMM} was the Buchberger-M\"oller's paper (cf. \cite{BM}). Authors of \cite{Rob} used the approach of \cite{BM} and some ideas of \cite{FGLM} for an efficient algorithm to zero-dimensional schemes in both affine and
projective spaces.  
In \cite{BBMo} similar ideas of using a generalized FGLM algorithm as a pattern algorithm were presented in order to compute Gr\"obner basis of ideals of free finitely generated algebras. In particular, it is introduced the pattern algorithm for monoid and group algebras
% using linear algebra in order to compute Gr\"obner bases or to have a Gr\"obner bases %framework for solving related problems 
. In \cite{alon-mora} a more general pattern algorithm which works on modules is introduced, many things behind of this idea of using linear algebra are formalized, notions like \lq\lq Gr\"obner technology" and \lq\lq Gr\"obner representations" are used. There are other approches which also generalized similar ideas to some settings, behind of all these works is the essential fact of using linear algebra techniques to compute in \lq\lq Gr\"obner bases schemes". 

The main goal of this paper is to show the application of techniques to linear codes like the ones in FGLM and subsequent works, which comes from an specification of the pattern algorithm for monoid algebras given in \cite{BBMo}, i.e. by taking an algebra associated to a linear code.

\section{Preliminaries}\label{pre}
The case of the algebra associated to a linear code we are going to introduce is connected with an ideal of a free commutative algebra; therefore, we will restric ourselves to the formulation of a pattern algorithm for a free commutative algebra. Similar settings can be performed in a free associated algebra or over modules (see \cite{BBMo,alon-mora}). 

Let $X:=\{x_{1},\dots,x_{n}\}$ be a finite set of variables, $[ X]$ the free commutative monoid on $X$, $K$ a field, $I$ an ideal of $K[ X]$, $I(F)$ the ideal of $K[ X]$ generated by $F\subset K[ X]$, $K[ X]/I$ the residue class algebra of $K[ X]$ module $I$. Let us denote by $1$ the empty word in $[ X]$, $L(u)$ the length of the word $u\in [X]$, and $Card(C)$ the cardinal of the set $C$. Let now $\prec$ be a semigroup total well ordering on $[ X ]$
(such an ordering is also called admissible), then for $f \in K [ X ]
\setminus \{0\}$, $T_{\prec}(f)$ is the maximal term of $f$ with respect
to $\prec$, $LC_{\prec}(f)$ is the leading coefficient of $f$ with respect to
$\prec$. Similarly, for $F \subset K [ X ]$, $T_\prec\{F\}$ is
the set of maximal terms of non-zero polynomials in $F$, $T_\prec(F)$ is
the semigroup ideal generated by $T\{F\}$. Moreover, for
the sake of simplicity in notation, $U_\prec(F)$ will be used instead of
$U(T_\prec(F))$, where $U$ lies in $\{G, N, B, I\}$. Of course, 
given an ideal $I$ and two different admissible orderings $\prec_1$ and $\prec_2$, in
general we have $U(T_{\prec_1}(I)) \neq U(T_{\prec_2}(I))$.
Notwithstanding this strong dependency on $\prec$, while a single admissible 
ordering $\prec$ is considered, so that no confusion arise,
we will often simply write $U(F)$ for $U_{\prec}(F)$.

Let $\tau \subset [ X ]$  be a semigroup
ideal of $[ X ]$, i.e., for $u \in [ X ]$
and $t \in \tau$, $tu \in \tau$. Then, it is well known that $\tau$
has a unique subset $G(\tau)$ of irredundant generators (probably
infinite). In the case of $I$ a zero-dimensional ideal, for $\tau=T(I)$, $G(\tau)$ is always finite. We are going to introduce for $\tau$ some notation and
terminology, which are similar to those introduced in \cite{MMM}.\\[0,1cm]
%\centerline{$Pred(w) : = \{u \in [X]\,\mid\, \exists\; x\in X\,
%(w=ux)\} \mbox{ (The set of predecessors of $w$)}.$}
%\vfill
\noindent
\hspace*{0.2cm}\begin{tabular}{ll}
%&\\
$Pred(w) : = \{u \in [X]\,\mid\, \exists\; x\in X\,
(w=ux)\}$ &$\,$ (the set of {\em predecessors of $w$}),\\
$N(\tau) := \{s \in [ X ] \mid s \notin
\tau\}$ &$\,$ ({\em outside of $\tau$}),\\
$B(\tau) := \{w \in \tau \mid Pred(w) \subset N(\tau)\}$ &$\,$
({\em border of $\tau$}),\\
$I(\tau) := \tau \setminus B(\tau)$ &$\,$
({\em interior of $\tau$}). \end{tabular} 

We remark that $w \in \tau$ lies in $G(\tau)$ if and only
if all its proper divisors are in $N(\tau)$ (that is if $Pred(w)\subset N(\tau)$). In the following
proposition, some basic results concerning $\tau$ and its regions are
summarized. Although they are very easy to prove, their importance is
crucial for FGLM techniques. 
\begin{proposition}[Properties of the semigroup ideal regions]\label{p:ireg} 
\begin{enumerate}
\item[{\it i.} ]For each $w \in \tau$ there exist $u \in [
X ]$ and $v \in B(\tau)$ s.t. $w = vu$.
\item[{\it ii.} ]For $x \in X$:
\begin{enumerate}
\item If $u \in N(\tau)$, then $ux \in N(\tau) \cup
B(\tau)$.
\item If $u \in B(\tau)$, then $ux \in
B(\tau) \cup I(\tau)$.   
\item If $u \in I(\tau)$, then $ux \in I(\tau)$.
\end{enumerate}
\item[{\it iii.} ]$N(\tau), N(\tau) \cup G(\tau), N(\tau) \cup
B(\tau)$ are order ideals, i.e., if $u$ belongs to one of these
subsets and $v$ divides $u$, then $v$ also belongs to the corresponding sets.
\end{enumerate} 
\end{proposition}

%\addtocounter{theo}{1}
\begin{theorem}[The vector space of canonical forms modulo an ideal]\label{t:vsp-cf} Let
$Span_{K}(N_\prec(I))$ be the $K$-vector space whose basis is $N_\prec(I)$. Then the
following holds: \begin{enumerate}
\item[{\it i.} ]$K \langle X \rangle = I \oplus Span_{K}(N_\prec(I))$ (this sum
is considered as a direct sum of vector spaces).
\item[{\it ii.} ]For each $f \in K [ X ]$ there is a unique
polynomial of $Span_{K}(N_\prec(I))$, denoted by $Can(f, I, \prec)$
such that $f - Can(f, I, \prec) \in I$; moreover:
\begin{enumerate}
\item $Can(f, I, \prec) = Can(g, I, \prec)$ if and only if $f - g \in
I$.
\item $Can(f, I, \prec) = 0$ if and only if $f \in I$.
\end{enumerate}
\item[{\it iii.} ]There is a $K$-vector space isomorphism between
$K [ X ]/I$ and $Span_{K}(N_\prec(I))$ (the isomorphism
associates the class of $f$ modulo $I$ with the canonical form $Can(f, I, \prec)$
of $f$ modulo $I$).
\end{enumerate} \end{theorem}
$Can(f, I, \prec)$ is called the canonical form of $f$
modulo $I$. We use simply $Can(f, I)$ if the ordering used is clear from the context.

We assume the readers to be familiar with definition and properties of Gr\"obner bases (see \cite{adams} for an easy to read introduction to Gr\"obner bases).

\begin{proposition}[Characterization of zero-dimensional ideals]\label{p:chzid} Let $G$ be a Gr\"obner basis of $I$ with respect to $\prec$. Then, $I$ is a zero-dimensional ideal (i.e. $dim_{K} K
[ X ]/I<\infty$) if and only if $N_\prec(G)$ is finite. Moreover, in such
a case, $dim_{K} K [ X ]/I = Card(N_\prec(G))$.
\end{proposition} 

\begin{definition}[Border basis] The border
basis of $I$ with respect to $\prec$ is the subset ${\cal B}(I,\prec) \subset I$ defined by:  

\centerline{ 
${\cal B}(I,\prec):=\{w-Can(w,I,\prec)\mid w\in B_{\prec}(I)\}
\mbox{\hspace{0.2cm} (the ${\cal B}$-basis of $I$)}.$}
\end{definition} 

Note that the ${\cal B}$-basis of $I$ is a Gr\"obner basis of $I$ that contains the reduced Gr\"obner basis.
\subsection{Matphi matrices and Gr\"obner representation}\label{sb:mat-Groeb}
The word Matphi appears by the first time in \cite{FGLM} to denote a procedure that computes a set of matrices (called matphi matrices) s.t. there is one matrix for each variable in $X$ and they describe the multiplication structure of the quotient algebra $K[X]/I$, where $I$ is a zero dimensional ideal. We often refer to this set of matrices as the matphi structure.

\begin{definition}[Gr\"obner representation, Matphi structure]\label{d:math-strc}
Let $I$ be a zero-dimensional ideal of $K[X]$, let $s=dim(K[X]/I)$. A Gr\"obner representation of $I$ is a pair $(N,\phi)$ consisting of
\begin{enumerate}
\item[{\it i. }] $N=\{N_1,\ldots,N_s\}$ s.t. $K[X]/I=Span_K(N)$, and
\item[{\it ii. }] $\phi:=\{\phi(k)\,\mid\, 1\leq k\leq n\}$, where $\phi(k)$ are the square matrices $\phi(k):=(a_{ij}^k)_{ij}$ s.t. for all $1\leq i\leq s$, $N_ix_k\equiv_I \sum_j\, a_{ij}^kN_j$.
\end{enumerate}
$\phi$ is called the matphi structure and the $\phi(k)$'s the matphi matrices. 
\end{definition}

See \cite{alon-mora} for a more general treatment of these concepts. Note that the matphi structure is indepent of the particular set $N$ of representative elements of the quotient $K[X]/I$. In addition, the matphi matrices allow to obtain the class of any product of the form $N_ix_k$ as a combination of the representative elements (i.e. as a linear combination of the basis $N$ for the vector space $K[X]/I$). 

\section{The FGLM pattern algorithm}\label{s:fglm_alg}
In this section we present a generalization of the FGLM algorithm for free commutative algebras, which allows to solve many different problems and not only the clasic FGLM convertion problem. The procedure we are presenting is
based on a sort of black box pattern: in fact, the description of the
steps 5 and 6 is only made in terms of their input and output. 
More precisely, we are assuming that  a term ordering ${\prec_1}$ 
is fixed on $[ X ]$, $I$ is a zero-dimensional
ideal (without this restriction the algorithm does not terminate), and
that the $K$-vector space $Span_{K}(N_{\prec_1}(I))$ is represented by
giving
\begin{itemize} 
\item[$\bullet$]a $K$-vector space $E$ which is endowed of
an {\em effective} function 
\[\mbox{{\bf LinearDependency}}[v,\{v_{1},\ldots,
v_{r}\}]\] which, for each finite set $\{v_{1},\ldots,v_{r}\} \subset
E$ of linearly independent vectors and for each vector $v\in E$,
returns the value defined by 
\[\left\{\begin{array}{ll} 
\{\lambda_{1},\ldots,\lambda_{r}\}\subset K & 
\mbox{if\ } v=\sum_{i=1}^{r}\lambda_{i}v_{i},\\
\mbox{\bf False} & \mbox{if $v$ is not a linear combination of
$\{v_{1},\ldots,v_{r}\}.$} \end{array} \right.\] 
\item[$\bullet$] an injective morphism $\xi : Span_{K}(N_{\prec_1}(I)) \mapsto E$. 
\end{itemize}

This informal approach allows a free choice of a suitable
representation of the space $Span_{K}(N_{\prec_1}(I))$ regarding an efficient
implementation of these techniques and a better complexity. Moreover, as an aside effect, it
enables us to present this generalization in such a way that it can
be applied on several more particular patterns and helps to make key
ideas behind the FGLM algorithm easier to understand. Let us start making some references to some subroutines of the
algorithm.

\noindent 
${\bf InsertNexts}[w, List, \prec]$ inserts properly the products $wx$
(for $x \in X$) in $List$ and sorts it by increasing ordering with
respect
to the ordering $\prec$. The reader should remark that ${\bf InsertNexts}$ could count the number of times that an element $w$ is inserted in $List$, so $w \in N_\prec(I) \cup T_\prec\{G\}$ if and only if this number coincide with the number of variables in the support of $w$, otherwise, it means that $w\in T_\prec(I)\setminus T_\prec\{G\}$, see \cite{FGLM}, this criteria can be used to know the boolean value of the test condition in Step~4 of the Algorithm~\ref{a:alg1}.

\noindent
${\bf NextTerm}[List]$ removes the first element from $List$ and
returns it.
%\pagebreak[4]
\begin{algorithm}[FGLM pattern algorithm]\label{a:alg1}\noindent$\,$\\
{\bf Input:} \hspace*{0.05cm} 
\begin{minipage}[t]{0.87\textwidth}
$\prec_2$, a term ordering on $[ X ]$; $\xi:Span_{K}(N_{\prec_1}(I)) 
\mapsto E$. 
\end{minipage}\\ 
%\nopagebreak[4]
\noindent
{\bf Output:} $rGb(I, \prec_2)$, the reduced Gr\"obner basis of $I$ w.r.t. the ordering $\prec_2$.\\[0.1cm] 
{\bf 1.} $G:=\emptyset$; $List:=\{1\}$; $N:=\emptyset$; $r:=0$;\\
{\bf 2.} {\bf While} $List \neq \emptyset$ {\bf do}\\
{\bf 3.} \hspace*{0.1cm} $w:={\bf  NextTerm}[List]$;\\
{\bf 4.} \hspace*{0.1cm} {\bf If} $w \notin T_{\prec_2}(G)\;$ (if $w$ is not a multiple of any element in $G$)
{\bf then} 
%\hspace{0.07cm}({\em it occurs iff }
%$t=1\mbox{ or }lr(t)\in N$)\footnote{Since the function InsertNexts produces only
%elements in $N_<(I)\cup rB_<(I)$ and because of Theorem 3.ii.}
\\
{\bf 5.} \hspace*{0.5cm} $v:=\xi(Can(w,I,\prec_1))$;\\
{\bf 6.} \hspace*{0.5cm} $\Lambda:={\bf LinearDependency}[v,
\{v_{1},\ldots,v_{r}\}]$;\\
{\bf 7.}\hspace*{0.6cm} $\mbox{\bf If }\Lambda\neq\mbox{ False }\mbox{ \bf
then }G:=G\cup\{w-\sum_{i=1}^{r}\lambda_{i}w_{i}\}$ \hspace{0.05cm}
({\em where } $\Lambda=(\lambda_1,\ldots,\lambda_r)$)\\
{\bf 8.} \hspace*{3.08cm} {\bf else }$\;r:=r+1$;\\
{\bf 9.} \hspace*{4cm} $v_{r}:=v$;\\
{\bf 10.} \hspace*{3.8cm} $w_{r}:=w; \; N:=N\cup\{w_{r}\}$;\\
{\bf 11.} \hspace*{3.7cm} $List:=\mbox{\bf InsertNexts}[w_{r},
List, \prec_2]$; \\{\bf 12.} $\mbox{\bf Return}[G]$.
\end{algorithm}

\begin{remark}\label{r:alg1}
\begin{enumerate} 
\item[{\it i. }]A key idea in algorithms like FGLM is to use the
relationship between membership to an ideal $I$ and linear
dependency modulo $I$, namely $\forall\; c_i\in K, s_i\in K[
X]$:

\noindent
\centerline{
$\sum_{i=1}^{r}c_{i}s_{i} \in I \setminus \{0\} \Longleftrightarrow
\{s_{1},\ldots,s_{r}\}\mbox{ is linearly dependent modulo $I.$}$}

\noindent
This connection with linear algebra was used for the firts time in
Gr\"obner bases theory since the very begining (see  \cite{B}).
\item[{\it ii. }] Since each element of $N_{\prec_2}(I)\cup B_{\prec_2}(I)$ belongs to $List$ at some moments of the algorithm and $List \subset N_{\prec_2}(I)\cup B_{\prec_2}(I)$ at each iteration of the algorithm, it is clear that one can compute ${\cal B}(I,{\prec_2})$ or the Gr\"obner representation $(N_{\prec_2}(I),\phi)$ of $I$ just by eliminating Step~4 of the algorithm and doing from Step~5 to Step~11 with very little changes in order to built those structures instead of $rGb(I,{\prec_2})$. 
\item[{\it iii. }]Note that Step~5 and 6 depends on the particular setting. In Step~5 it is necessary to have a way of computing $Can(w,I,\prec_1)$ and the corresponding element in $E$, while in Step~6 we need an effective method to decide linear dependency.
\item[{\it iv. }]Complexity analysis of this pattern algorithm can be found in \cite{BBMo} for the more general case of free associative algebras, and for a more general setting in \cite{alon-mora}. Of course, having a pattern algorithm as a model, it is expected that for particular applications, one could do modification and specification of the steps in order to improve the speed and decrease the complexity of the algorithm by taking advantage of the particular structures involved.
\end{enumerate}
\end{remark}
\subsection{The change of orderings: a particular case}\label{sb:fglm_part}
Suppose we have an initial ordering $\prec_1$ and the reduced Gr\"obner basis of $I$ for this ordering, now we want to compute by the FGLM algorithm the new reduced Gr\"obner basis for a new ordering $\prec_2$. Then the vector space $E$ is $K^s$, where $s=dim(K[X]/I)$. In Step~5, $Can(w,I,\prec_1)$ can be computed using the reduced Gr\"obner basis $rGb(I,\prec_1)$ and the coefficients of this canonical form build the vector of $E$ corresponding to this element (the image by the morphism $\xi$). Then Step~6 is perfomed using pure linear algebra. 
\section{FGLM algorithm for monoid rings}\label{s:mon}
The pattern algorithm is presented in \cite{BBMo} for the free monoid algebra, we will restrict here to the commutative case. Let $M$ be a finite commutative monoid generated by
$g_{1},\ldots,g_{n}$; $\xi:[ X ] \rightarrow M$,
the canonical morphism that sends $x_{i}$ to $g_{i}$; $\sigma \subset 
[ X ]\times[ X ]$, a presentation of $M$
defined by $\xi$  ($\sigma=\{(w,v)\,\mid \, \xi(w)=\xi(v)\}$).
Then, it is known that the monoid ring $K[M]$ is isomorphic to $K
[ X ]/I(\sigma)$, where $I(\sigma)$ is the ideal generated by $P(\sigma)=\{w-v\, \mid\,(w,v)\in\sigma\}$; moreover, any
Gr\"obner basis $G$ of $I(\sigma)$ is also formed by binomials of
the above form. In addition, it can be proved that $\{(w,v)\,\mid\,w-v\in
G\}$ is another presentation of $M$.

Note that $M$ is finite if and only if $I=I(\sigma)$ is zero-dimensional. We will show that in order to compute $rGb(I)$, the border basis or the Gr\"obner representation of $I$, one only needs to have $M$ given by a concrete
representation that allows the user to multiply words on its
generators; for instance: $M$ may be given by permutations, matrices
over a finite field, or by a more abstract way (a complete or
convergent presentation). Accordingly, we are going to do the
necessary modifications on Algorithm~\ref{a:alg1} for this case. 

We should remark that in this case $\prec_1=\prec_2$, then at the begining of the algorithm the set $N_{\prec_1}(I)$ is unkown (which is not the case of the change of orderings). It could be precisely a goal of the algorithm to compute a set of representative elements for the quotient algebra. 

Now consider the natural extension of $\xi$ to an algebra morphism ($\xi:K[X]\mapsto K[M]$), note that the restriction of $\xi$ to $Span_{K}(N_{\prec_1}(I))$ ( $\xi:Span_{K}(N_{\prec_1}(I))\mapsto K[M]$) is an injective morphism; moreover,  $\xi(w)=\xi(Can(w,I,\prec_1))$, for all $w \in [X]$. Therefore, the image of $Can(w,I,\prec_1)$ can be computed as $\xi(w)$, and the linear dependency checking will find out whether $w$ is a new canonical form (i.e. $w \in N_{\prec_1}(I)$) or not (i.e. $w\in T_{\prec_1}(rGb(I,\prec_1))$). Hence, Step 5 will be

\centerline{
$v:=\xi(u)g_{i},\mbox{ where } u\in Pred(w) \mbox{ and } ux_{i}=w.$}

Moreover, let $w_1,\ldots,w_r$ be elements of $N_{\prec_1}(I)$ and $v_i=\xi(w_i)$, for $1\leq i\leq r$. Then ${\mbox{\bf LinearDependency}}[v,\{v_{1},\ldots,v_{r}\}]$
can be computed as

\centerline{
$\left\{\begin{array}{ll} 
%\xi^{-1}(v)
v_j  & 
\mbox{ if\ } v=v_j, \mbox{ for some } j\in[1,r],
% }\in \{v_{1},\ldots,v_{r}\}
\\
\mbox{\bf False} & \mbox{ otherwise}.  \end{array} \right.$}

Finally, Step 7 changes into: 

\centerline{
$\mbox{\bf If } \Lambda \neq \mbox{ \bf  False } {\bf\ then\ } G:=G\cup
\{w-w_j\}.$}

\begin{remark}\label{r:mon_case}
\begin{enumerate} 
\item[{\it i. }]This example shows that the capability of the
$K$-vector space $E$ w.r.t. {\bf LinearDependency}, that is demanded
in the Algorithm 1, is required only on those sets of vectors
$\{v_{1},\ldots, v_{r},v\}$ that are built in the algorithm, which means in this case that {\bf LinearDependency} is reduced to the {\bf Member} checking, i.e., $v$ is linear dependent of $\{v_1,\ldots,v_r\}$ if and only if it belongs to this set. 
%\item[{\it ii. }]We can use $\xi^{-1}$, because for any element $v\in List$ we know %$\xi(v)$, and for the elements $v_i$, $1\leq i\leq r$ we know its corresponding element %$w_i$ ($v_i=\xi(w_i)$).
\item[{\it ii. }]When a word $w$ is analyzed by the algorithm, all the elements in $Pred(w)$ have been already analyzed ($\xi(u)$ is known for any $u\in Pred(w)$), this is the case whenever $\prec_1$ is an admissible ordering. Therefore, the computation of $\xi(w)$ is immediate.
\end{enumerate} 
\end{remark}

We will show the case of linear codes as a concrete  setting for an application of the FGLM pattern algorithm for monoid rings, where the monoid is given by a set of generators and a way of multiply them.
\section{FGLM algorithm for linear codes}\label{s:fglm-bin}
For the sake of simplicity we will stay in the case of binary linear codes, where more powerfull structures for applications are obtainned as an output of the corresponding FGLM algorithm (for a general setting see \cite{bbw-rep,bbm}). From now on we will refer to linear codes simply as codes.

Let $\mathbb{F}_2$ be  the finite field with $2$ elements. Let $\mathcal C$ be a binary code of dimension $k$ and length $n$ ($k\leq n$), so that the $n\times (n-k)$ matrix $H$ is a {\em parity check matrix} ($c\cdot H=0\,$ if and only if $c\in\mathcal C$). Let $d$ be the minimum distance of the code, and $t$ the error-correcting capability of the code 
($t=\left[\frac{d-1}{2}\right]$, where $[x]$ denotes the greater integer less than $x$). Let $B({\mathcal C},t)=\{y\in \mathbb{F}_2^n\mid \exists\; c\in\mathcal{C}\;(d(c,y)\leq t) \}$, it 
is well known that the equation $eH=yH$ 
has a unique solution $e$ with $\mathrm{weight}(e)\leq t$, for $y\in B({\mathcal C},t)$. 

%\subsection{The monoid associated with a binary code}\label{sb:mon-bc}
Let us consider the free commutative monoid $[X]$ 
generated by the $n$ variables $X:=\{ x_{1},\dots ,  x_{n}\}$. We have the following map from $X$ to $ 
\mathbb{F}_2^n$:
%\begin{equation}\begin{split}
$\psi: X\to \mathbb{F}_2^n$, where $x_{i}\mapsto\, e_i$ (the $i$-th coordinate vector). 
%\end{split}\end{equation}
The map $\psi$ can be extended in a natural way to a morphism from  $[X]$ onto  
$\mathbb{F}_2^n$, where 
$\psi( \prod_{i=1}^n\, x_{i}^{\beta_{i}})=(\beta_{1}\,\mathrm{mod}\;  2 ,\dots,
\beta_{n}\,\mathrm{mod}\;  2 )$.
%\end{equation}

A binary code $\mathcal C$ defines an equivalence relation $R_{\mathcal C}$ in  
$\mathbb{F}_2^n$ given by 
%\begin{equation}\label{eq:eq-rel}
$(x,y)\in R_{\mathcal C}$ if and only if $x-y\in \mathcal C$.
%\end{equation}
If we define $\xi (u):=\psi (u)H$, where $u\in [X]$, the above 
congruence can be translated to $[X]$ by the morphism $\psi$ as  
%\begin{equation}
$u\equiv_{\mathcal C} w$ if and only if $(\psi(u),\psi(w))\in R_{\mathcal 
C}$, that is, if $\xi(u)=\xi(w)$. The morphism $\xi$ represents the 
transition of the syndromes from $\mathbb{F}_2^n$ to $[X]$; therefore, $\xi(w)$ 
is the \lq\lq syndrome" of $w$, which is equal to the syndrome of $\psi(w)$.

\begin{definition}[The ideal associated with a binary code]\label{d:id-cod} Let $\mathcal C$ be a binary code. 
The ideal $I({\mathcal C})$  associated with $\mathcal C$ is
%\begin{equation} 

\centerline{$
I(\mathcal C):=\langle \{w -u \,\mid\, \xi(w)=\xi(u)
\}
\rangle\subset K[X].
$}
%\end{equation}
\end{definition}

\section{The algorithm for binary codes}\label{sb:fglm-code}
The monoid $M$ is set to be $\mathbb{F}_2^{n-k}$ (where the syndromes belong to). Doing $g_i:=\xi(x_i)$, note that $M=\mathbb{F}_2^{n-k}=\langle g_1,\ldots, g_n \rangle$. Moreover, $\sigma:=R_{\mathcal C}$, hence $I(\sigma)=I({\mathcal C})$.  Let $\prec$ be an admissible ordering. Then the FGLM algorithm for linear codes can be used to compute the reduced Gr\"obner basis, the border basis, or the Gr\"obner representation for $\prec$.

\begin{algorithm}[FGLM for binary codes]\label{a:reducida}$\quad$\newline 
{\bf Input:} $n,H$ the parameters for a given binary code, $\prec$ an admissible ordering.\newline 
{\bf Output:} $rGb(I({\mathcal C}),\prec)$.
\begin{itemize} 
\item[{\bf 1.}] $List:=\{1\}, N:=\emptyset ,r:=0,G=\{\}$; 
\item[{\bf 2.}] {\bf While} $List\neq \emptyset$ do 
\item[{\bf 3.}]  $\qquad w:=\mathbf{NextTerm}[List]$; 
\item[{\bf 4.}] $\quad\;\;$  If $w \notin T(G)$; 
\item[{\bf 5.}]  $\qquad v:=\xi (w)$; 
\item[{\bf 6.}]  $\qquad \Lambda:= \mathbf{Member} [v ,\{v_1,\dots ,v_r\}]$; 
\item[{\bf 7.}] $\qquad${\bf If} $\Lambda \neq$ False $\;${\bf then} $G:=G \cup \{w-w_j\}$; 
%\item[8.] $\qquad$ $\qquad\qquad \phi(u,x_k):=w_j$ 
\item[{\bf 8.}] $\qquad$ {\bf else} $r:=r+1$; 
\item[{\bf 9.}] $\qquad\qquad\, v_r:=v^\prime$; 
\item[{\bf 10.}]  $\qquad\qquad\, w_r:=w,\; N:=N\cup \{w_r\}$;  
\item[{\bf 11.}]  $\qquad\qquad List:=\mathbf{InsertNext}[w_r,List]$;  
\item[{\bf 12.}] {\bf Return}$[G]$.
%\end{center}  
\end{itemize} 
\end{algorithm}

In many cases of FGLM applications a good choice of the ordering $\prec$ is a crucial point in order to solve a particular problem. In the following theorem it is shown the importance of using a total degree compatible ordering (for example the Degree Reverse Lexicographic). Let us denote by $<_T$ a total degree compatible ordering.

\begin{theorem}[Canonical forms of the vectors in $B(C,t)$]\label{t:DecGB}
Let $\mathcal C$ be a code and let $G_T$ be the reduced Gr\"ober
basis with respect to $<_T$. If $w \in [X]$ satisfies
$\,\mathrm{weight}(\psi(Can(w,G_T))) \leq t$ then
$\psi(Can(w,G_T))$ is the error vector corresponding to
$\psi(w)$. On the other hand, if
$\mathrm{weight}(\psi(Can(w,G_T)))
> t$ then $\psi(w)$ contains more than $t$ errors.
\end{theorem}
\begin{proof}
If we assume that $\mathrm{weight}(\psi(Can(w,G_T))) \leq t$ then, we can infer at once that $\psi(w)\in B({\mathcal 
C},t)$ and $\psi(Can(w,G_T))$ is its error vector (notice that $\xi(w)=\xi(Can(w,G_T))$ and the unicity of the error vector). 

Now, if $\mathrm{weight}(\psi(Can(w,G_T))) > t$, we have to prove that $\psi(w) \notin B({\mathcal 
C},t)$. It is equivalent to show that $\mathrm{weight}(\psi(Can(w,G_T)))~\leq~t\,$ if  $\,\psi(w) \in B({ \mathcal 
C},t)$. Let $\psi(w)$ be an element of $B({\mathcal 
C},t)$ and let $e$ be its error vector then, $\mathrm{weight}(e) \leq t$. Let $\mathrm{w}_e$ be the squarefree representation of $e$. Note that  $\mathrm{weight}(e)$ coincides with the total degree of $\mathrm{w}_e$; accordingly, $L(\mathrm{w}_e) \leq t$. On the other hand, $Can(w,G_T)<_T \mathrm{w}_e$, which implies that $L(Can(w,G_T))\leq L(\mathrm{w}_e)$ (because $<_T$ is degree compatible). Hence, $\mathrm{weight}(\psi(Can(w,G_T)))\leq L(Can(w,G_T))\leq t$. $\Box$
\end{proof}

The computation of the error-correcting cability of the code $t$ can be done in the computing process of Algorithm~\ref{a:reducida}  (see the example in Section~\ref{e:fglm-cod} and \cite{bbw-rep}). The previous theorem allows us to use the computed reduced Gr\"obner basis for solving the decoding problem in general binary codes, but also with such a powerful tool available, it is expected to be able to study the structure of the codes, like some combinatorics properties. Some possible examples are the permutation-equivalence of codes (see \cite{bbm}), and some problems related with binary codes associated with the set of cycles in a graph (finding the set of minimal cycles and a minimal cycle basis of the cycles of a graph), see \cite{bbmf}.

To generalize Theorem~\ref{t:DecGB} for non binary linear codes have some conflicts with the needed ordering; however, the FGLM algorithm can be still used to compute the border basis or a Gr\"obner representation for the ideal $I({\mathcal C})$ and it will be possible to solve the problems that one can solve with the reduced Gr\"obner basis in the case of binary codes. Those problems are explained in \cite{bbw-rep}. In addition, \cite{bbm} contains some results and examples about the application of this setting to general linear codes and, in binary codes, for studying the problems of decoding and the permutation-equivalence. 
\subsection{An example}\label{e:fglm-cod}
Let ${\mathcal C}$ be the linear code over $\mathbb{F}_2^6$ determined by 
the parity check matrix $H$ given below. 
The set ${\mathcal C}$ of codewords is given in the right hand side. 
The minimum distance is $d=3$, so, $t=1$, the numbers of variables is $6$, 
$<_T$ is set to be the Degree Reverse Lexicographic ordering with $x_{i+1} >_T x_i$. 
Only essential parts of the computation will be described.

\begin{tabular}{ll}
$H=\left|\begin{array}{ccc}
  1\hspace{0.3 cm} & 1\hspace{0.3 cm} & 1\hspace{0.3 cm}\\
  1\hspace{0.3 cm} & 0\hspace{0.3 cm} & 1\hspace{0.3 cm}\\
  0\hspace{0.3 cm} & 1\hspace{0.3 cm} & 1\hspace{0.3 cm}\\
  1\hspace{0.3 cm} & 0\hspace{0.3 cm} & 0\hspace{0.3 cm}\\
  0\hspace{0.3 cm} & 1\hspace{0.3 cm} & 0\hspace{0.3 cm}\\
  0\hspace{0.3 cm} & 0\hspace{0.3 cm} & 1\hspace{0.3 cm}\\
  \end{array} \right|$&,
\begin{minipage}[t!]{0.6\textwidth}
$\begin{array}{l}
C=\{(0,\,0,\,0,\,0,\,0,\,0),\,(1,\,0,\,1,\,1,\,0,\, 0),
    \\ \hspace*{0.2 cm} 
    (1,\,1,\,1,\,0,\,0,\,1),\,(1,\,1,\,0,\,0,\,1,\,0),
    (0,\,1,\,0,\,1,\,0,\,1),\\ \hspace*{0.2 cm} 
    (0,\,1,\,1,\,1,\,1,\,0),\,(0,\,0,\,1,\,0,\,1,\,1),\,(1,\,0,\,0,\,1,\,1,\,1) \}.\\
\end{array}$
\end{minipage}
\end{tabular}

%\newpage
\noindent{\bf Application of Algorithm \ref{a:reducida}: }$List:=\{1\}$; $N:=\{\}$; $r:=0$; $w:=1$; $\xi(1)=(0,0,0)$; 
$N:=N \cup\{1\}=\{1\}$;
$\xi(N):=\{(0,0,0)\}$; $List:=\{x_{1},x_{2},x_{3},x_{4},x_{5},x_{6} \}$; 
$w:=x_{1}$; 
$\xi(x_{1})=(1,1,1)$; $N:=\{1,x_{1}\}$; $\xi(N):=\{(0,0,0),(1,1,1)\}$;\\
After analyzing $x_6$ we are at the following stage:\\
$N:=\{1,x_1,x_2,x_3,x_4,x_5,x_6\},$   and $List=\{x_1^2, x_1x_2, x_1x_3, x_1x_4, x_1x_5, x_1x_6, x_2^2,$ $x_2x_3, x_2x_4, x_2x_5, x_2x_6, x_3^2, x_3x_4, x_3x_5, x_3x_6, x_4^2, x_4x_5, x_4x_6, x_5^2, x_5x_6, x_6^2\}$.

There is still one element left in $N$ because there are 7 elements in $N$ 
of a total of 8 ($2^{6-3}$). Taking the elements of $List$ from $x_1^2$ to $x_1x_5$ they are a linear combination of elements already in $N$ (their syndromes are in the list of syndromes computed $\xi(N)$). Therefore, $G:=\{x_1^2-1,x_1x_2-x_5,x_1x_3-x_4,x_1x_4-x_3,x_1x_5-x_2\}$, for example $x_1x_2-x_5$ is obtained, bacause when $w=x_1x_2$, first note that $Pred(w) \subset N$, which means that it is either a new irreducible element or a head of a binomial of the reduced basis. Then $\xi(x_1x_2)$ is computed and we got that $\xi(x_1x_2)=\xi(x_5)$. This means that $x_1x_2-x_5$ belongs to $G$. Also $x_1x_2$ is the first minimal representation which is not in $N$, this implies that $t=\mathrm{weight}(\psi(x_1x_2))-1$ (see \cite{bbw-rep}). The next element in $List$, $w=x_1x_6$, is the last element that will be included in $N$ and the corresponding multiples will be included in $List$. 

From this point, the algorithm will just take elements from $List$ and it analyzes in each case whether it is in $T\{rGb(I({\mathcal C}),<_T)\}$ (like $x_2x_3$) or in\linebreak[4] $T(rGb(I({\mathcal C}),<_T)\setminus T\{rGb(I({\mathcal C}),<_T)\}$ (like $x_1x_2x_6$), this process is executed until the $List$ is empty when the last element $x_1x_6^2$ of the list is analyzed. Finally, the reduced Gr\"obner basis for $<_T$ is 

$G:=\{x_1^2-1,x_1x_2-x_5,x_1x_3-x_4,x_1x_4-x_3,x_1x_5-x_2,
x_2^2-1, x_2x_3-x_1x_6, x_2x_4-x_6, x_2x_5-x_1, x_2x_6-x_4, x_3^2-1, x_3x_4-x_1, x_3x_5-x_6, x_3x_6-x_5, x_4^2-1, x_4x_5-x_1x_6, x_6x_4-x_2, x_5^2-1, x_5x_6-x_3, x_6^2-1
\}.$

Now let us assume that a vector $y=(1,1,1,0,1,0)$ is received, the corresponding word is $w=x_1x_2x_3x_5$. Then we compute $w_e=Can(w,G)=x_3$. As $\mathrm{weight}(\psi(w_e))=1$ ($t=1$ is the error-correcting capability); therefore, the error vector is $e=(0,0,1,0,0,0)$, and the codeword is $c=(1,1,0,0,1,0)$.

\end{document}